\documentclass{amsart}

\usepackage{amsmath,amsthm,amssymb,amscd,enumerate,mathtools,enumitem,tikz-cd,bm}
\usepackage{amsfonts}
\usepackage{rotating}
\usepackage{euscript}
\usepackage{pst-node}
\usepackage{epsfig,verbatim}
\def\be{\begin{equation}}
\def\ee{\end{equation}}

\def\C{{\mathbb C}} 
\def\f{\EuScript}
 
\def\P{{\mathbb P}}

\def\R{{\mathbb R}} 
\def\Q{{\mathbb Q}}

\def\phi{{\varphi}}

\def\deg{{\rm deg\,}}

\def\Aut{{\rm Aut}}

\def\bp{\begin{proposition}}
\def\ep{\end{proposition}}

\def\bt{\begin{theorem}}
\def\et{\end{theorem}}
\def\br{\begin{remark}}
\def\er{\end{remark}}
\def\be{\begin{equation}}
\def\bee{\begin{equation*}}
\def\l{\label}

\def\ee{\end{equation}}
\def\eee{\end{equation*}}
\def\bl{\begin{lemma}}
\def\el{\end{lemma}}
\def\bc{\begin{corollary}}
\def\ec{\end{corollary}}
\def\pr{\noindent{\it Proof. }}

\def\bd{\begin{definition}}
\def\ed{\end{definition}}
\def\t{\widetilde}
\def\h{\widehat}
\def\hat{\widehat}

\def\t{\widetilde }

\newtheorem{theorem}{Theorem}[section]
\newtheorem{lemma}[theorem]{Lemma}
\newtheorem{definition}[theorem]{Definition}
\newtheorem{corollary}[theorem]{Corollary}
\newtheorem{proposition}[theorem]{Proposition}
\newtheorem{problem}[theorem]{Problem}

\theoremstyle{definition}

\theoremstyle{definition}
\newtheorem{remark}[theorem]{Remark}
\def\bpr{\begin{problem}}
\def\epr{\end{problem}}

\mathtoolsset{showonlyrefs}

\begin{document}

\title[]{Holomorphic maps sharing preimages over finitely generated fields
}

\author[F. Pakovich]{Fedor Pakovich}
\thanks{
This research was supported by ISF Grant No. 1092/22} 
\address{Department of Mathematics, 
Ben Gurion University of the Negev, P.O.B. 653, Beer Sheva,  8410501, Israel}
\email{pakovich@math.bgu.ac.il}

\begin{abstract} 
Let \( R \) be a compact Riemann surface, and let  
\( P: R \to \P^1(\C) \) and  
\( Q: R \to \P^1(\C) \) be holomorphic maps. In this paper, we investigate the following problem: under what conditions do the preimages \( P^{-1}(K) \) and \( Q^{-1}(K) \) coincide for some infinite set $K$ contained in $\P^1({\bm k})$, where $\bm k$ is a finitely generated subfield of $\C$  (e.g., a number field)?
Equivalently, we study holomorphic correspondences that admit an infinite completely invariant set contained in $\P^1({\bm k})$. We show that if such a set exists, then there is a  
holomorphic Galois covering \( \Theta: R_0 \to \P^1(\C) \),  
where \( R_0 \) has genus zero or one, such that \( P \) and \( Q \) are ``compositional left factors" of \( \Theta \). We also consider a more general   equation  \( P^{-1}(K_1) = Q^{-1}(K_2) \), where $K_1$ and $K_2$ are infinite subsets of $\P^1({\bm k})$.  
 \end{abstract}


\maketitle

\section{Introduction} 
The famous five-point-theorem of Nevanlinna \cite{nev} states that if $f$ and $g$ are  two non-constant meromorphic  functions on $\C$ such that the preimages $f^{-1}(\{a\})$ and  $g^{-1}(\{a\})$ coincide for five distinct points $a$ of $\P^1(\C),$ then $f\equiv g.$  Nevanlinna obtained this result as an application of 
the deep theory of the distribution of values of meromorphic functions that he developed. On the other hand, if we assume that $f$ and $g$ are {\it rational} functions on $\P^1(\C)$, then an easy application of the Riemann-Hurwitz formula implies that 
$f\equiv g$ whenever the above equalities hold for any  four distinct points of $\P^1(\C)$ (see \cite{adam}, \cite{piz}, \cite{sau}).  

The problem of describing rational functions that share preimages becomes much more subtle when one considers the preimages of \textit{sets} rather than those of individual points. This difficulty arises already in the case of polynomials. 
For example, the following problem, posed in \cite{y}, remained open for nearly twenty years. Let  
$P$ and $Q$ be non-constant polynomials of the same degree such that  
$$P^{-1}(\{-1,1\})=Q^{-1}(\{-1,1\}).$$ Does it follow that 
$P=\pm Q$? This problem was affirmatively solved in \cite{p0}. Further results 
regarding polynomials satisfying the condition 
 \be \l{00} P^{-1}(K)=Q^{-1}(K)\ee for some {\it compact} set  $K$  in $\C$ were obtained in 
\cite{d}, \cite{d2}, \cite{p1}, \cite{p2}, \cite{p3}. 

A description of solutions to a more general  equation  
\be \l{0} P^{-1}(K_1) = Q^{-1}(K_2),\ee  
where \( P \) and \( Q \) are non-constant polynomials and \( K_1 \) and \( K_2 \) are arbitrary compact subsets of \( \mathbb{C} \), not necessarily equal, was obtained in \cite{p3}. Specifically, in \cite{p3}, condition \eqref{0} was related to the functional equation  
\begin{equation}  
\label{2}  
F \circ P = G \circ Q,  
\end{equation}  
in polynomials. It is clear that for any polynomial solution of \eqref{2} and any compact set \( K \subset \mathbb{C} \), one obtains a solution of \eqref{0} by setting  
\begin{equation}  
\label{3}  
K_1 = F^{-1}(K), \quad K_2 = G^{-1}(K).  
\end{equation}  
Somewhat unexpectedly, the main result of \cite{p3} states that {\it all} solutions to \eqref{0} can be constructed in this way, provided that the compact set defined by either side of \eqref{0} contains at least \( \text{LCM}(\deg P, \deg Q) \) points. This condition, in particular, holds for any infinite \( K_1 \) and \( K_2 \). Since Ritt’s theory of polynomial decompositions \cite{r1} provides a quite precise description of solutions to \eqref{2}, we thus obtain a precise description of solutions to \eqref{0}. This description can then be applied to various related problems, including the classification of polynomials that share a Julia set, the description of commuting and semiconjugate polynomials, and the study of invariant curves for polynomial endomorphisms of \(\P^1(\C) \times \P^1(\C)\) (see \cite{p3}, \cite{pj}, \cite{int} for further details).

The methods employed in the aforementioned papers are restricted to the polynomial case, and the problem of finding solutions to \eqref{00} and \eqref{0} when $P$ and $Q$ are arbitrary \textit{rational} functions and $K_1$ and $K_2$ are compact subsets of $\C$ or $\P^1(\C)$ remains largely unresolved. With the exception of the note \cite{des}, where the solutions of the  equation  
$$
\beta_1^{-1}[-1,1] = \beta_2^{-1}[-1,1]
$$  
were described in the case where $\beta_1$ and $\beta_2$ are rational Belyi functions, 
the only related paper known to us is \cite{bel}, where the broader problem of describing invariant sets of correspondences was investigated.

We recall that a \textit{correspondence}, defined by an algebraic curve 
\be \l{cur}
f(x, y) = 0,
\ee
over $\C$, is a multivalued map that assigns to a point $x$ the set of all points $y_i$ satisfying $f(x, y_i) = 0$.  
Along with this forward map, one also defines a backward map that assigns to a point $y$ the set of all points $x_j$ such that $f(x_j, y) = 0$.   
If $R$ is a desingularization of the curve \eqref{cur}, and $z \mapsto (P(z), Q(z))$ is its parametrization given by a pair of holomorphic maps on $R$, then the forward and backward maps of the correspondence defined by \eqref{cur} are described by the multivalued functions $Q(P^{-1}(z))$ and $P(Q^{-1}(z))$, respectively. Thus,   a pair of holomorphic maps $P$ and $Q$ can be  regarded as the primary object of investigation.

The correspondences can be naturally composed, and their dynamics have been extensively studied from various perspectives (see, for instance, the recent works \cite{bl}, \cite{d3}, \cite{in}, \cite{lm}, and the references therein). In this context, describing the sets satisfying \eqref{00} for a pair of holomorphic maps $P$ and $Q$ on a compact Riemann surface $R$ is clearly equivalent to describing the completely invariant sets of the associated correspondences.   In particular, the case where the maps 
$P$ and $Q$ are rational functions corresponds to the situation in which the curve \eqref{cur} has genus zero.


When applied to equation \eqref{00}, the main result of \cite{bel} states that if \( P:R\rightarrow \P^1(\C) \) and \( Q:R\rightarrow\P^1(\C) \) are holomorphic maps  on a compact Riemann surface $R$ satisfying \eqref{00} for {\it infinitely many finite sets} \( K \), then there exists a rational function \( F \) such that  
\be \label{4}  
F \circ P = F \circ Q.  
\ee  
Observe that if \eqref{4} holds, then \eqref{00} is satisfied for every set of the form \( K = F^{-1}(\hat{K}) \), where \( \hat{K} \subset \P^1(\C) \), and, in particular, for every fiber of \( F \). Consequently, all irreducible completely invariant sets of the associated correspondence are finite.
Notice also the  \eqref{4} obviously implies that \be \l{rav} \deg P=\deg Q.\ee 

In this paper, we study equation \eqref{00}, assuming that \( P: R \to \P^1(\C) \) and \linebreak \( Q: R \to \P^1(\C) \) are holomorphic maps  on a compact Riemann surface $R$, but replacing the condition that \( K\) is compact with the assumption that \( K \) is an infinite set contained in $\P^1({\bm k})$, where $\bm k$ is a finitely generated subfield of $\C$  (e.g., a number field). In other words, we investigate correspondences that have an infinite, completely invariant set contained in $\P^1({\bm k})$. We also 
study equation \eqref{0} under the assumption that $K_1$ and $K_2$ are infinite sets contained in $\P^1({\bm k})$. 
Notice that \eqref{0} can also be interpreted in terms of correspondences: it expresses 
the condition that \(K_1\) is mapped to \(K_2\) by the forward map, while \(K_2\) is 
mapped to \(K_1\) by the backward map of the correspondence associated with the pair \(P, Q\).

We remark that describing solutions of  \eqref{0}, where \( K_1 \) and \( K_2 \) are subsets of \( \mathbb{C} \) satisfying certain restrictions, reduces to describing solutions for which the maps \( P \) and \( Q \) have no non-trivial common compositional right factor in the following sense: the equalities 
\be \l{lu}  
P = \t{P} \circ W, \quad Q = \t{Q} \circ W,   
\ee  
where 
\be
\t{P}: \t{R} \to \P^1(\C), \quad \t{Q}: \t{R} \to \P^1(\C),\ \text{ and }\  W: R \to \t{R}, 
\ee
are holomorphic maps between compact Riemann surfaces, imply that \( \deg W = 1 \). Indeed, it follows from \eqref{0} and \eqref{lu} that
\[
W^{-1}(\t{P}^{-1}(K_1)) = W^{-1}(\t{Q}^{-1}(K_2)),
\]
which implies
\[
\t{P}^{-1}(K_1) = \t{Q}^{-1}(K_2).
\]
Thus, any solution $P$, $Q$ of \eqref{0} reduces to a solution   \( \t{P} \), \( \t{Q} \), where \( \t{P} \) and \( \t{Q} \) have no non-trivial common compositional right factor, and we will primarily focus on such solutions.

In connection with the problem under consideration, we mention the recent note \cite{hei}, where it was shown, using properties of height functions, that if $P$ and $Q$ are non-constant rational functions over $\C$, and $K$ is an infinite subset of $\mathbb{P}^1({\bm k})$, where $\bm{k}$ a finitely generated subfield of $\C$, then the condition  
\[
P^{-1}(K) \subseteq Q^{-1}(K)
\]  
implies  
\[
\deg Q \ge \deg P.
\]  
In particular, equality \eqref{00} implies equality \eqref{rav}. Our approach is different; however, it does not immediately yield \eqref{rav}, even in the case when the holomorphic maps $P$ and $Q$ are rational functions on $\C\P^1$.  
Thus, the result of \cite{hei} is complementary to ours.

Our main result concerning equation \eqref{0} is the following statement.

\bt \l{t0} Let \( R \) be a compact Riemann surface, and let \( P: R \to \P^1(\C) \) and \( Q: R \to \P^1(\C) \) be holomorphic maps having no non-trivial common compositional right factor. Assume that the equality
\[
P^{-1}(K_1) = Q^{-1}(K_2)
\]
holds for some infinite sets \( K_1 \) and \( K_2 \) contained in $\P^1({\bm k})$, where $\bm k$ is a finitely generated subfield of $\C$. 
 Then, there exist compact Riemann surfaces \( R_1 \) and \( R_2 \) of genus at most one and holomorphic Galois coverings \( \Theta_1: R_1 \to \P^1(\C) \) and \( \Theta_2: R_2 \to \P^1(\C) \) such that the equalities
\[
\Theta_1 = P \circ U, \quad \Theta_2 = Q \circ V
\]
hold for some holomorphic maps \( U: R_1 \to R \) and \( V: R_2 \to R \).

\et 

Theorem \ref{t0} implies, in particular, that \( g(R) \leq 1\). More importantly, since Galois coverings \( \Theta: R \to \P^1(\C) \) with \( g(R) \leq 1 \) are easy to describe, Theorem \ref{t0} imposes strict constraints on \( P \) and \( Q\). In particular, they can exhibit only very limited branching (see Lemma \ref{21} below).

The simplest example of holomorphic maps \( P \) and \( Q \) satisfying the conclusion of Theorem \ref{t1} is any pair of rational Galois coverings. Moreover, we show that for any such pair, one can construct sets \( K_1 \) and \( K_2 \) contained in  $\P^1({\bm k})$, where $\bm k$ is a finitely generated subfield of $\C$, such that \eqref{0} holds. 
These examples demonstrate the existence of pairs \( P \), \( Q \) for which \eqref{0} holds, but \eqref{2} does not. Indeed, denoting the groups of covering transformations of \(P\) and \(Q\) by \(\Gamma_P\) and \(\Gamma_Q\), we see that if \eqref{2} holds, then the rational function defined by any part of this equality must be invariant under the group \(\Gamma = \langle \Gamma_P, \Gamma_Q \rangle\). This is possible only if \(\Gamma\) is finite, which is not necessarily the case. 
On the other hand, if \( \Gamma \) {\it is} finite, then such solutions are indeed obtained from the convenient formulas \eqref{2} and \eqref{3} (see Section \ref{exm}).  


Our second main result is the following specification of Theorem \ref{t1} in the case $K_1=K_2.$

\bt \l{t1} Let \( R \) be a compact Riemann surface, and let \( P: R \to \P^1(\C) \) and \( Q: R \to \P^1(\C) \) be holomorphic maps having no non-trivial common compositional right factor. Assume that the equality
\[
P^{-1}(K) = Q^{-1}(K)
\]
holds for an infinite set \( K \)  contained in $\P^1({\bm k})$, where $\bm k$ is a finitely generated subfield of $\C$.   Then, there exist a compact Riemann surface \( R_0 \) of genus at most one and a holomorphic Galois covering \( \Theta: R_0 \to \P^1(\C) \) such that the equalities
\[
\Theta = P \circ U, \quad \Theta = Q \circ V
\]
hold for some holomorphic maps \( U: R_0 \to R \) and \( V: R_0 \to R \).
\et 
For Theorem \ref{t1}, the simplest examples of maps \(P\) and \(Q\) that satisfy its conclusion are any rational Galois coverings \(P\) and \(Q\) such that
\[
Q = P \circ \mu
\]  
for some \(\mu \in \Aut(\P^1(\C))\). Moreover, for such \(P\) and \(Q\), one can construct a set \(K\) contained in a finitely generated field such that \eqref{00} holds, and these examples generally do not reduce to \eqref{4} (see Section \ref{exm}).

For illustration, we consider the following simple example. Let
\[
P = z^2, \quad Q = (z+1)^2 = P \circ (z+1),  
\]  
and let \(K\) be the set of squares of integers. Then we  have
\[
P^{-1}(K) = Q^{-1}(K) = \mathbb{Z}.  
\]  
However, the equality \eqref{4} is impossible, because, if it held, the rational function defined by any part of this equality would be invariant under the transformation \(z \mapsto z+1\).

This paper is organized as follows. In the second section, we begin by recalling several definitions and results related to fiber products and normalizations, and then prove two results concerning equation \eqref{00}, extending Theorem \ref{t0}. In the third section, we prove several results related to equation \eqref{0} and establish Theorem \ref{t1}, again in an extended form.  Finally, in the fourth section, we construct examples of solutions to \eqref{00} and \eqref{0}, illustrating Theorems \ref{t0} and \ref{t1}.

\section{Holomorphic maps sharing preimages of different sets}

\subsection{Fiber products and normalizations} 
In this subsection, we review several definitions and results concerning fiber products and normalizations.

Let  $V_i: E \rightarrow R_i$,  $1\leq i \leq k,$ 
where $k\geq 2,$ be holomorphic maps between  compact Riemann surfaces. We say that the maps  $V_i$, $1 \leq i \leq k,$
 have  no non-trivial common compositional right factor if the 
equalities 
\be \l{ho} V_i=\t V_i\circ W, \ \ \  1 \leq i \leq k,\ee where $W:   E\rightarrow \t E$ and 
$\t V_i: {\t E}\rightarrow R_i,$ $1 \leq i \leq k,$ are holomorphic  maps 
between compact Riemann surfaces 
imply that $\deg W=1.$ Denoting by  $\f M(S)$ the field of meromorphic functions on a compact Riemann surface $S$,  this condition 
can be restated 
as the requirement 
\be \l{des} \f M( E)=V_1^*(\f M({R}_1))\cdot V_2^*(\f M({R}_2))\cdot\ \dots \ \cdot V_k^*(\f M( {R}_k)),\ee 
meaning that the field $\f M( R)$ is the compositum of its subfields $V_i^*(\f M(E_i)),$ \linebreak $1 \leq i \leq k.$

Let us recall that if  $P_i: R_i \rightarrow \P^1(\C)$, $1 \leq i \leq k,$  are holomorphic maps between compact Riemann surfaces, then the fiber product of $P_i$, $1 \leq i \leq k,$ is 
 a  collection
\be \l{nota} P_1 \times P_2 \times \dots \times P_k=\bigcup\limits_{j=1}^{n}\{ E_j,V_{j1},V_{j2}, \dots, V_{jk}\},\ee 
where $n=n(P_1,P_2,\dots ,P_k)$ is an integer positive number and $ E_j,$ $1\leq j \leq n,$ are compact Riemann surfaces provided with holomorphic maps
$$V_{ji}:\,  E_j\rightarrow  R_i, \ \ \  1\leq i \leq k, $$
such that  \be \l{dee} P_1\circ V_{j1}=P_2\circ V_{j2}=\dots = P_k\circ V_{jk}, \ \ \ 1\leq j \leq n,\ee 
and for any holomorphic maps $T_i:\, { E}\rightarrow  R_i,$  $1\leq i \leq k, $
between compact Riemann surfaces satisfying \be \l{bua}  P_1\circ T_1=P_2\circ T_2=\dots = P_k\circ T_k\ee  there exist a uniquely defined  index $j$, $1\leq j \leq n$, and 
a holomorphic map \linebreak $W:\, { E}\rightarrow  E_j$ such that
\be \l{ae} T_i= V_{ji}\circ  W, \ \ \ 1\leq i \leq k.\ee  

Notice that the definition implies that for every $j,$ $1\leq j \leq n,$ the maps \linebreak $V_{ji}:\,  E_j\rightarrow  R,$ $1\leq i \leq n$,  have no non-trivial common compositional  right factor. 
In the other direction, if $T_i:\, { E}\rightarrow  R_i,$  $1\leq i \leq k, $ are holomorphic maps 
between compact Riemann surfaces satisfying \eqref{bua}  and having no non-trivial common compositional  right factor, then   \eqref{ae} holds for some $j$, $1\leq j \leq n$, and isomorphism $W:\,   E\rightarrow  E_j.$ 
We will call each collection $\{E_j, V_{j1}, V_{j2}, \dots, V_{jk}\}$, $1 \leq j \leq n$, a component of the fiber product $P_1 \times P_2 \times \dots \times P_k.$ We call the genus of a component  the genus of $E_j$, $1 \leq j \leq n$.

The fiber product  is defined in a unique way up to natural isomorphisms, and 
can be described by the following algebro-geometric construction. Let us consider the algebraic variety 
\be \l{ccuurr} L =\{(x_1,x_2,\dots,x_k)\in  R_1\times  R_2\times \dots \times  R_k \, \vert \,  P_1(x_1)=P_k(x_2)=\dots =P_k(x_k)\}.\ee
Let us denote by $L _j,$ $1\leq j \leq n$,  irreducible components of $L $, by 
$ E_j$,  $1\leq j \leq n$, their desingularizations, 
 and by $$\pi_j:  E_j\rightarrow L _j, \ \ \ 1\leq j \leq n,$$ the desingularization maps.
Then the compositions  $$x_i\circ \pi_j:  E_j\rightarrow  R_i, \ \ 1\leq i \leq k,\ \  1\leq j \leq n,$$ 
extend to holomorphic maps
$$V_j:\,  E_j\rightarrow  R_i, \ \ 1\leq i \leq k,\ \  1\leq j \leq n,$$
and the collection $\bigcup\limits_{j=1}^{n}\{ E_j,V_{j1},V_{j2}, \dots, V_{jk}\}$ is   the fiber product $P_1 \times P_2 \times \dots \times P_k.$

Let $R$ be a compact  Riemann surface and 
$Q:\, R\rightarrow \P^1(\C)$  a holomorphic map. 
A  {\it normalization} of $Q$ is defined as a compact Riemann surface  $N_Q$ together with a holomorphic Galois covering  of the 
minimal degree $\t Q:N_Q\rightarrow \P^1(\C)$ such that \linebreak
 $\t Q=Q\circ H$ for some  holomorphic map $H:\,N_Q\rightarrow R$. 
The normalization is characterized by the property that  the field extension 
$$\f M(N_Q)/\t Q^*(\f M(\P^1(\C)))$$ is isomorphic to the Galois closure 
of the extension $$\f M(R)/Q^*(\f M(\P^1(\C))).$$

  We  recall that  {\it an orbifold} $\f O$ on $\P^1(\C)$ is a ramification function $\nu:\P^1(\C)\rightarrow \mathbb N$, 
which takes the value $\nu(z)=1$ except at finite points, and  
 {\it the   Euler characteristic} of $\f O$ is the number
$$ \chi(\f O)=2+\sum_{z\in \P^1(\C)}\left(\frac{1}{\nu(z)}-1\right).$$  
{\it The signature} \( \nu(\f O) \) of \( \f O \) is defined as the list of all values \( \nu(z) \), where \( z \) ranges over the points of \( \f O \) with \( \nu(z) > 1 \), and each value is included as many times as it occurs.
With each holomorphic map \( Q:\, R\rightarrow \P^1(\C) \) between compact Riemann surfaces, one can associate its {\it ramification orbifold} $\f O^Q$  by setting  \( \nu^Q(z) \) equal to the least common multiple of the local degrees of \( Q \) at the points of the preimage \( Q^{-1}\{z\} \). 

The following statement provides different characterizations of the condition appeared in Theorem \ref{t0} and Theorem \ref{t1}.

\bl \l{21} Let $R$  be a compact  Riemann surface, and $P:R\rightarrow \P^1(\C)$ a holomorphic map of degree at least two. Then the following conditions are equivalent.
\begin{enumerate}[label=(\roman*)]  
\item
There exist a compact Riemann surfaces $R_0$ of genus at most one and a holomorphic Galois covering \( \Theta: R_0 \to \P^1(\C) \)  such that $\Theta = P \circ V$ for some  holomorphic map \( V: R_0 \to R \). 
\item 
The inequality  $g(N_P)\leq 1$ holds. 
\item 
The inequality $\chi(\f O^P)\geq 0$  holds. 
\item The signature $\nu(\f O^P)$ belongs to the list 
\be \l{list}\{2,2,2,2\} \ \ \ \{3,3,3\}, \ \ \  \{2,4,4\}, \ \ \  \{2,3,6\}, \ee 
 \be \l{list2} \{l,l\}, \ \ l\geq 2,  \ \ \ \{2,2,l\}, \ \ l\geq 2,  \ \ \ \{2,3,3\}, \ \ \ \{2,3,4\}, \ \ \ \{2,3,5\}.\ee
\end{enumerate}  
\el
\pr 
The equivalency $i)\Leftrightarrow ii)$ follows from the minimality of the Galois covering $\t Q$ and the fact that the genus does not increase under holomorphic maps.  
The equivalency $ii)\Leftrightarrow iii)$ follows from the Riemann-Hurwitz formula 
(see e.g.  \cite{plo}, Lemma 3.1). Finally, the equivalency $iii)\Leftrightarrow iv)$ is  obtained by a direct calculation (see e.g. \cite{fk},  
IV.9.3, IV.9.5). \qed 

 The normalization  of a holomorphic map \( P:R\rightarrow \P^1(\C) \) of degree \( d \)  
can be described in terms of the fiber product of \( P \) with itself \( d \) times as follows (see \cite{fried}, $\S$I.G, or \cite{plo}, Section 2.2).  
Consider the algebraic variety  
\[
L^P =\{(x_1,x_2,\dots,x_k)\in  R\times  R\times \dots \times  R \, \vert \,  P(x_1)=P(x_2)=\dots =P(x_k)\}.
\]  
Let \( \hat{ L}^{P} \) be a variety obtained from \( L^{P} \) by removing the components contained in varieties of the form  
\[
x_{j_1}=x_{j_2},   \quad  1\leq j_1,j_2 \leq d, \quad j_1\neq j_2,
\]  
and let \( N \) be an irreducible component of \( \hat{ L}^{P} \).  
Further, let \( \pi^{\prime}: {N_P}\rightarrow N \) be the desingularization map, and  
\( \t P: {N_P}\rightarrow \P^1(\C) \) a holomorphic map induced by the composition  
\[
{N_P}\xrightarrow{\pi^{\prime}} {N}\xrightarrow{\pi_i}\P^1(\C)\xrightarrow{X}\P^1(\C),
\]
where $\pi_i$ is the projection to any coordinate.

 In this notation, the following statement holds. 
 
\bt \l{fr} The map $\t P: {N_P}\rightarrow  R$  is the normalization of $P$. \qed 
\et


\subsection{Proof of Theorem \ref{t0} and its extensions} 
We deduce Theorem \ref{t0} from the following result, which  
follows from the Faltings theorem. 

\bt \l{mm}
Let $R$ be a compact Riemann surface, and let  \( P_i: R \to \P^1(\C) \), \( 1 \leq i \leq k \), where \( k \geq 2 \), be holomorphic maps having no non-trivial common compositional right factor.  
Assume that for some infinite set \( K \) contained in $\P^1({\bm k})$, where $\bm k$ is a finitely generated subfield of $\C$, the inclusions  
\begin{equation} \label{int} 
P_i(P_k^{-1}(K)) \subset \P^1({\bm k}), \quad 1 \leq i \leq k-1, 
\end{equation}  
hold. Then \( g(N_{P_k}) \leq 1 \). 
\et
\pr 
The image of $ R$ under the map \be \l{the} \theta:
z\rightarrow (P_1(z),P_2(z),\dots ,P_k(z))\ee  
defines an irreducible algebraic curve $X $ in the space $(\P^1(\C))^k$ with coordinates 
$(x_1,x_2,\dots ,x_k)$.  Moreover, since $P_i$, $1\leq  i\leq k,$  have no non-trivial common compositional right factor, $ R$ is the desingularization of $X$. 

Setting \( d = \deg P_k \), let us consider the product  
\( X^d \) in the space \( (\P^1(\C))^{kd} \) with coordinates  
\[
(x^1_1, x^1_2, \dots, x^1_k, x^2_1, x^2_2, \dots, x^2_k, \dots, x^d_1, x^d_2, \dots, x^d_k)
\]
and define the algebraic variety \( E \) as the intersection of \( X^d \) and the algebraic variety  
\[
x^1_k = x^2_k = \dots = x^d_k
\]
in \( (\P^1(\C))^{kd} \).  
By adjoining, if necessary, finitely many coefficients of the equations defining \( E \)  
in \( (\P^1(\C))^{kd} \) to \( \bm{k} \), we may assume that \( E \) is defined over \( \bm{k} \).

Since the map \eqref{the} is an isomorphism off a finite set, the map 
\be \l{und} (\theta, \theta, \dots, \theta):R^d\rightarrow X^d\subset (\P^1(\C)^{kd}\ee  
induces an isomorphism between components of the curve $\hat{ L}^{P_k}\subset R^d$ defined above and components of the curve $ E\subset X^n$.   On the other hand,  it is easy to see that 
the  condition of the theorem implies that the image of $\hat{ L}^{P_k}$ under \eqref{und} 
 has infinitely many points over   $\bm k$. Therefore, at least one of irreducible components of this image also has infinitely many points over $\bm k$. Now the statement    of the theorem  follows from Theorem \ref{fr} and the Faltings theorem, which states that 
if an irreducible algebraic curve $C$ defined over a finitely generated field $\bm k$ of characteristic zero has infinitely many $\bm{k}$-points, 
then $g(C)\leq 1$   (\cite{fa}). \qed

The proof of Theorem \ref{mm} given above is a modification of the proof of one of the implications of the main result in \cite{alg}, which characterizes algebraic curves \( X : F(x, y) = 0 \) defined over \(\overline{\mathbb{Q}}\) that satisfy the following property: there exist a number field \( \bm{k} \) and an infinite set \( S \subset \bm{k} \) such that, for every \( y \in S \), the roots of the polynomial \( F(x, y) \) belong to \( \bm{k} \). The connection between these problems becomes evident upon noting that if a curve $X$ as above is parametrized by holomorphic maps \( \phi: R \to \mathbb{P}^1(\mathbb{C}) \) and \( \psi: R \to \mathbb{P}^1(\mathbb{C}) \), 
then \( \phi(\psi^{-1}(S)) \subset \bm{k} \).

Theorem \ref{mm} implies the following result, which extends Theorem \ref{t0}.

\bt \l{mmc} Let $R$ be a compact Riemann surface, and let \( P_i: R \to \P^1(\C) \), \( 1 \leq i \leq k \), where \( k \geq 2 \), be  holomorphic maps having no non-trivial common compositional right factor.   
 Assume that  the  
 equality 
\be\l{cr} P_1^{-1}(K_1)=P_2^{-1}(K_2)=\dots =P_k^{-1}(K_k)\ee  holds for some 
infinite sets $K_i$, $1\leq i \leq l,$ contained in $\P^1({\bm k})$, where $\bm k$ is a finitely generated subfield of $\C$.   Then $g(N_{P_i})\leq 1,$ $1\leq i \leq k.$
\et 
\pr Since equality \eqref{cr} implies that for every \( j \), \( 1 \leq j \leq k \), the inclusion  
\[
P_i(P_j^{-1}(K_j)) \subset \P^1({\bm k}), \quad 1 \leq i \leq k, \quad i \neq j,  
\]  
holds, the statement of the theorem follows from Theorem \ref{mm}. \qed

\section{Holomorphic maps sharing preimages of the same set} 
In this section, using fiber products, we derive refined versions of Theorem \ref{mmc} under the additional assumption that the equality  
\be \label{u}  
K_1 = K_2 = \dots = K_k  
\ee  
holds in \eqref{cr}.

\bl\l{l3} Let $R$ be a compact Riemann surface, and let  
\( P_i: R \to \P^1(\C) \), \( 1 \leq i \leq k \), where \( k \geq 2 \), be holomorphic maps having no non-trivial common compositional right factor.  
Assume that the equality  
\be \label{it0}  
K' = P_1^{-1}(K) = P_2^{-1}(K) = \dots = P_k^{-1}(K)  
\ee  
holds for some subsets \( K \subset \P^1(\C) \) and \( K' \subset R \). 
Then  for any holomorphic maps  $V_i:\, E\rightarrow  R,$  $1\leq i \leq k, $
between compact Riemann surfaces satisfying \be \l{bua2}  P_1\circ V_1=P_2\circ V_2=\dots = P_k\circ V_k\ee
the equalities 
\be \l{it00} V_1^{-1}(K')=V_2^{-1}(K')=\dots =V_k^{-1}(K')\ee 
and 
\be \l{f} (P_{i_1}\circ V_{j_1})^{-1}(K)=(P_{i_2}\circ V_{j_2})^{-1}(K), \ \ \ 1\leq i_1,i_2,j_1,j_2\leq k,\ee
 hold. 
\el
\pr 
Setting 
$$F= P_1\circ V_1=P_2\circ V_2=\dots = P_k\circ V_k,$$ we see that 
$$V_i^{-1}(K')=(P_i\circ V_i)^{-1}(K)=F^{-1}(K), \ \ \ 1\leq i \leq k.$$ 
Equality \eqref{f} follows now from \eqref{it0} and \eqref{it00}. \qed

\bl \l{l}
Let $ R_1$, $ R_2$, $ R_3$  be compact Riemann surfaces, $V_i: R_1 \rightarrow  R_2$, \linebreak $1\leq i \leq n,$  holomorphic maps having  no non-trivial common compositional right factor, and   $U_j: R_2 \rightarrow  R_3$, $1\leq j \leq m,$ 
 holomorphic maps having  no non-trivial common compositional right factor.  Then 
the holomorphic maps $U_j\circ V_i:  R_1 \rightarrow  R_3$, $1\leq i \leq n,$ $1\leq j \leq m,$ also have  no non-trivial common compositional right factor. 
\el
\pr 
Let $K$ be the compositum of the fields  $$(U_j\circ V_i)^*(\f M( R_3)),\ \ \  1\leq i \leq n, \ \ 1\leq j \leq m,$$ and  $K_j$, $1\leq j \leq k,$   the compositum
$$K_{j}=(U_1\circ V_j)^*(\f M( R_3))\cdot (U_2\circ V_j)^*(\f M( R_3))\cdot\ \dots \ \cdot (U_m\circ V_j)^*(\f M( R_3)), \ \ \ 1\leq j \leq n.$$
By the condition, 
$$U_1^*(\f M( R_3))\cdot U_2^*(\f M( R_3))\cdot\ \dots \ \cdot U_m^*(\f M( R_3))=\f M( R_2),$$ 
implying that
$$ K_j=V_j^*(\f M( R_2)), \ \ \ 1\leq j \leq n.$$ 
Since $K$ contains $ K_j,$ $1\leq j \leq n,$ and 
$$V_1^*(\f M( R_2))\cdot V_2^*(\f M( R_2))\cdot\ \dots \ \cdot V_n^*(\f M( R_2))=\f M( R_3)$$ by the condition, this implies that $K=\f M( R_3).$ \qed

Since for any holomorphic maps \( P: R \to \P^1(\C) \) and \( V: \t R \to R \) between 
compact Riemann surfaces the inequality 
\( g(N_P) \leq g(N_{P \circ V}) \) holds, the following result may be viewed 
as a stronger version of Theorem \ref{mmc} in the case where \eqref{u} holds.

\bt\l{m3}  Let $R$ be a compact Riemann surface, and let \( P_i: R \to \P^1(\C) \), \( 1 \leq i \leq k \), where \( k \geq 2 \), be holomorphic maps having no non-trivial common 
compositional right factor.   Assume that the equality 
\be \l{it} P_1^{-1}(K)=P_2^{-1}(K)=\dots =P_k^{-1}(K)\ee  holds for some 
infinite set $K$ contained in $\P^1({\bm k})$, where $\bm k$ is a finitely generated subfield of $\C$. 
Then for    every  component $ \{E,V_1,V_2,\dots ,V_k\}$ of the fiber product  \linebreak  $P_1 \times P_2 \times \dots \times P_k$ 
the inequalities $$g(N_{P_i\circ V_j})\leq 1, \ \ \ 1\leq i,j\leq k,$$ hold. 
\et
\pr Since \(P_i\), \(1\leq i \leq k\), have no non-trivial common compositional right factor by assumption and the same holds for the maps \(V_i\), \(1\leq i \leq k\), as they form a component of a fiber product, it follows from Lemma \ref{l} that the holomorphic maps \(P_i \circ V_j\), \(1\leq i,j\leq k\), also have no non-trivial common compositional right factor. Moreover, by Lemma \ref{l3}, the equalities \eqref{bua2} imply the equalities \eqref{f}. Applying now Theorem \ref{mm} to the maps \(P_i \circ V_j\), \(1\leq i,j\leq k\), we conclude that \(g(N_{P_i \circ V_j})\leq 1\),  \(1\leq i,j\leq k\). \qed

For holomorphic maps $P_i,$ $1\leq i \leq k$, of degree $d_i$,  $1\leq i \leq k$,
let us consider a component $\{E,V_{i,j}, 1\leq i \leq k, 1\leq j \leq d_i\}$ of 
the fiber product 
\be \l{comp} \Pi=P_1 ^{\times d_1}\times P_2 ^{\times d_2}\times \dots \times P_k ^{\times d_k}\ee 
such that the corresponding irreducible component in the variety 
$$P_1(x_{11})=P_1(x_{12})=\dots =P_1(x_{1d_1})=P_2(x_{21})=\dots =P_2(x_{2d_2})=\dots $$
$$\dots = P_k(x_{k1})=\dots =P_k(x_{kd_k})$$ 
is not contained in a variety of the form $$x_{i,j_1}=x_{i,j_2}, \ \  1\leq i \leq k,  \ \  1\leq j_1,j_2 \leq d_i, \ \ j_1\neq j_2,$$ and set 
\be \l{po} F=P_i\circ V_{i,j}, \ \ 1\leq i \leq k, \ \ 1\leq j \leq d_i.\ee

\bt \l{else}  Let $R$ be a compact Riemann surface, and let \( P_i: R \to \P^1(\C) \), \( 1 \leq i \leq k \), where \( k \geq 2 \), be holomorphic maps having no non-trivial common 
compositional right factor.
Assume that the equality 
\be P_1^{-1}(K)=P_2^{-1}(K)=\dots =P_k^{-1}(K)\ee  holds for some 
infinite set $K$ contained in $\P^1({\bm k})$, where $\bm k$ is a finitely generated subfield of $\C$.  Then  $g(E)\leq 1$ and $F$ is a Galois covering.
 \et 
 \pr Since the functions $P_i$, $1\leq i \leq k,$ where each function is taken $d_i$ times still have no  non-trivial common compositional right factor, $g(E)$ is one or zero by Theorem \ref{m3}. Thus, we only must prove that  $F$ is a Galois covering. 
Let $\{\t E, U_1,U_2,\dots, U_k\}$ be a component of the fiber product $\t P_1\times \t P_2\times \dots \times \t P_k$, where  
$\t P_i$ is a normalization of $P_i$, $1\leq i \leq k,$ and  
\be \l{impl0} \t F=\t P_1\circ U_1=\t P_2\circ U_2=\dots =\t P_k\circ U_k. \ee 
Notice that since $\t P_i$, $1\leq i \leq k,$ are Galois coverings, all such components are isomorphic and $\t F$ is a Galois covering. Thus, to prove the theorem, it is enough to prove that $E$ and $\t E$ are isomorphic and the equality 
$F=\t F\circ \mu$ holds for some isomorphism $\mu:E\rightarrow \t E$. 

For each fixed \( i \), \( 1 \leq i \leq k \), the equalities \eqref{po} imply, by Theorem \ref{fr}, that  
\be \l{wou}  
V_{i,j} = \t{U}_{i,j} \circ W_i, \quad 1 \leq j \leq d_i,  
\ee  
where \( \{ \t{E}_i, \t{U}_{i,j}, 1 \leq j \leq d_i \} \) is a component of \( P_i^{\times d_i} \), and \( W_i : \t{E} \to \t{E}_i \) is a holomorphic map.  
Substituting \eqref{wou} into \eqref{po} for all \( i \), \( 1 \leq i \leq k \), we obtain  
\be \l{impl}  
F = \t{P}_1 \circ W_1 = \t{P}_2 \circ W_2 = \dots = \t{P}_d \circ W_d.  
\ee 
Moreover, \( W_1, W_2, \dots, W_k \) have no nontrivial common compositional right factor; otherwise, \eqref{wou} would imply that \( V_{i,j} \), \( 1\leq i \leq k \), \( 1\leq j \leq d_i \), have such a factor.  Thus, $\{E, W_1,W_2,\dots, W_k\}$ is a component of the fiber product $\t P_1\times \t P_2\times \dots \times \t P_k$ and  $
F = \t{F} \circ \mu  $ 
for some isomorphism \( \mu: E \to \t{E} \). \qed

Theorem \ref{else} implies the following statement extending Theorem \ref{t1} form the introduction. 

\bt   Let $R$ be a compact Riemann surface, and let 
$P_i: R \rightarrow \P^1(\C)$,  $1\leq i \leq k,$ 
where $k\geq 2,$ be holomorphic maps  having  no non-trivial common compositional right factor. 
Assume that the equality 
\be P_1^{-1}(K)=P_2^{-1}(K)=\dots =P_k^{-1}(K)\ee  holds for some 
infinite set $K$ contained in $\P^1({\bm k})$, where $\bm k$ is a finitely generated subfield of $\C$. 
 Then there exist a compact Riemann surface $R_0$  and a holomorphic Galois covering \( \Theta: R_0 \to \P^1(\C) \) such that the equalities
\[
\Theta = P_i \circ U_i,  \ \ \ 1\leq i \leq k,
\]
hold for some holomorphic maps \( U_i: R_0 \to R \),  $1\leq i \leq k.$
 \et 
\pr Since, by Theorem \ref{else}, the function \( F \) in \eqref{po} is a Galois covering, 
we can set \( \Theta = F \) and \( U_i = V_{i,1} \) for \( 1 \leq i \leq k \). \qed

\section{\l{exm} Examples of holomorphic maps sharing preimages} 
In this section, we construct examples of solutions to equations 
\eqref{00} and \eqref{0} illustrating Theorems \ref{t0} and \ref{t1}.

\bt \l{tt1} Let \( R \) be a compact Riemann surface of genus zero or one, 
and let 
$P_i: R \rightarrow \P^1(\C)$,  $1\leq i \leq k,$ 
where $k\geq 2,$ be holomorphic  Galois coverings. 
Then there exist  a finitely generated subfield $\bm k$ of $\C$ and 
infinite sets $K_i\subset \P^1({\bm k})$,  $1\leq i \leq k,$ 
such that the equality 
\be \l{esa} 
P_1^{-1}(K_1)=P_2^{-1}(K_2)=\dots =P_k^{-1}(K_k)
\ee 
holds.
\et 
\pr Assume first that $g(R)=0$, that is, $R=\P^1(\C)$. Let us  denote  by \( \Gamma_i \), $1\leq i \leq k,$  
the group of covering transformations of \( P_i \), $1\leq i \leq k,$   and consider 
the group 
$\Gamma$ generated by \( \Gamma_i \), $1\leq i \leq k,$ and some  M\"obius transformation of infinite order $\mu$, for example, $\mu(z)= z+1$.   Let us take now an arbitrary point $z_0\in \P^1(\C)$, and consider its orbit  
\be \l{bf} 
S = \bigcup_{\sigma \in \Gamma} \sigma(z_0)
\ee 
under the action of $\Gamma$. 
Notice that since $\mu\in \Gamma$, the set $S$ is infinite.

Since for every $z \in S$ the set $S$ contains the orbit of $z$ under 
the action of each group $\Gamma_i$, $1 \leq i \leq k$, there exist 
infinite sets $K_i \subset \P^1(\C)$, $1 \leq i \leq k$, such that 
$S = P_i^{-1}(K_i)$, $1 \leq i \leq k$.
 Moreover, it is clear that $S\subset\P^1({\bm k'})$, 
where $\bm k'$ is generated over $\Q$ by the coefficients of elements of \( \Gamma \),  and    $K_i\subset \P^1({\bm k})$, $1\leq i \leq k,$ where $\bm k$ is generated over $\bm k'$ by the coefficients of the rational functions \( P_i \), $1\leq i \leq k.$ Thus, $\bm k$ is finitely generated. Notice that since any finitely generated algebraic extension of $\Q$ is 
a number field, if $z_0$ and the coefficients of $P_i$, 
$1 \leq i \leq k$, are algebraic numbers, then the field $\bm{k}$ is a number field.

In case $g(R)=1$, the proof is modified as follows. Let us consider some holomopric maps $\phi$ and $\psi$ on $R$  having  no non-trivial common compositional right factor, and  a plane curve  $X:f(x,y)=0$ parametrized by  $\theta: z\rightarrow (\phi(z), \psi(z))$. Since $\phi$ and $\psi$  have  no non-trivial common compositional right factor, $R$ is the desingularization  of $X$, implying that for every holomorphic map 
$F : R \to \P^1(\C)$ there exists a rational function 
$\hat F(x,y) \in \C(x,y)$ such that 
$$
\hat F \circ \theta = F
$$
for all but finitely many points $z \in R$. 
Similarly, for every $\sigma \in \Aut(R)$ there exists a rational 
function $\hat \sigma(x,y)\in \C(x,y)$ such that 
$$
\hat \sigma \circ \theta = \theta \circ \sigma
$$
for all but finitely many points $z \in R$.

Let \( \Gamma \subset \Aut(R) \) be the group generated by the groups 
\( \Gamma_{P_i} \), \( 1 \leq i \leq k \), of covering transformations of \( P_i \), 
together with a shift \( \mu \) of infinite order on \( R \), and let 
\( \h{\Gamma} \subset \Aut(X) \) be the group consisting of all elements \( \h{\sigma} \), 
where \( \sigma \in \Gamma \). 
Let us take an arbitrary point  $(x_0,y_0)\in X(\C)$, and set  
\be\l{set}
\h S = \bigcup_{\h\sigma \in \h\Gamma} \h\sigma(x_0,y_0).
\ee
Since the set \( \h{S} \) is countable, we may, if necessary, change \( (x_0, y_0) \) 
so that the equalities  
\[
\h{P}_i \circ \theta = P_i, \quad 1 \leq i \leq k,
\]  
hold for every \( z \in S \). 
Thus, to prove the theorem it is enough to show that 
there exist  a finitely generated subfield $\bm k$ of $\C$ and 
infinite sets $K_i\subset \P^1({\bm k})$,  $1\leq i \leq k,$ 
such that the equality 
\[
\h P_1^{-1}(K_1)=\h P_2^{-1}(K_2)=\dots =\h P_k^{-1}(K_k)
\]
holds. Taking into account 
that the addition operation on $X$ is defined over the field of definition of $X$, this can be done 
 as in the first part of the proof. \qed
 \vskip 0.2cm
The formulation of Theorem \ref{tt1} does not require the maps  
\( P_i \), \( 1 \leq i \leq k \), to have no non-trivial common compositional right factor. 
On the other hand, since these maps are Galois coverings, they have such a factor 
if and only if the intersection  
\be \l{gr1}  
\Gamma_1 \cap \Gamma_2 \cap \dots \cap \Gamma_k  
\ee  
is nontrivial.
Furthermore, if  
\be \l{ho0}  
P_i = \widetilde{P}_i \circ W, \quad 1 \leq i \leq k,  
\ee  
where \( W: R \to \widetilde{R}_i \) and \( \widetilde{P}_i: \widetilde{R}_i \to E_i \),  
\( 1 \leq i \leq k \), are holomorphic maps between compact Riemann surfaces such that  
\( \deg W > 1 \) and the maps \( \widetilde{P}_i \), \( 1 \leq i \leq k \), have no nontrivial  
common compositional right factor, then the group \eqref{gr1} is  the group of covering transformations of \( W \).

Notice also that if the group  
\be \l{gr2}  
G = \langle \Gamma_1, \Gamma_2, \dots, \Gamma_k \rangle  
\ee  
is finite, then there exists a holomorphic Galois covering  
\( A: R \to \mathbb{P}^1(\mathbb{C}) \)  
such that \( G \) is its group of covering transformations.  
The inclusions \( \Gamma_i \subseteq G \), \( 1 \leq i \leq k \), then imply that  
\be \l{bn}  
A = F_1 \circ P_1 = F_2 \circ P_2 = \dots = F_k \circ P_k  
\ee  
for some rational functions \( F_i \), \( 1 \leq i \leq k \).  
Moreover, since \( S \) is a union of \( G \)-orbits, there exists a set \( K \) such that \( S = A^{-1}(K) \).  
Thus, 
$$  
S = P_i^{-1}(K_i) = P_i^{-1}(F_i^{-1}(K)), \quad 1 \leq i \leq k,  
$$  
implying  that  
\be \l{oof}  
K_i = F_i^{-1}(K), \quad 1 \leq i \leq k.  
\ee  

In the general case, for the solutions constructed in Theorem \ref{tt1}, 
rational functions \( F_i \), \( 1 \leq i \leq k \), satisfying \eqref{bn} do not exist, 
since \eqref{bn} would imply that \( G \) is finite, which need not hold.

\bt \l{tt2} Let \( R \) be a compact Riemann surface of genus zero or one, 
and let 
$P_i: R \rightarrow \P^1(\C)$,  $1\leq i \leq k,$ 
where $k\geq 2,$ be holomorphic  Galois coverings of the form 
$$P_i=P_k\circ \mu_i, \ \ \ 1\leq i \leq k-1,$$ where $\mu_i\in \Aut(R)$,   $1\leq i \leq k-1.$   
Then there exist  a finitely generated subfield $\bm k$ of $\C$ and 
an infinite set $K\subset \P^1({\bm k})$ such that the equality 
\[
P_1^{-1}(K)=P_2^{-1}(K)=\dots =P_k^{-1}(K)
\]
holds.
\et 
\pr 
The proof is obtained by a modification of the proof of Theorem \ref{tt1} as follows. Assuming first that $R=\P^1(\C)$, let us define a group $\Gamma\subset \Aut(\P^1(\C))$ as the group generated by the group $\Gamma_k$ of covering transformation of $P_k$, the M\"obius transformations 
$\mu_i$,   $1\leq i \leq k-1, $  and some  M\"obius transformation of infinite order $\mu$.  Let us take now a point $z_0\in \P^1(\C)$ and define a subset $S$ of $\P^1(\C)$ by the formula \eqref{bf}.   
Since the group \( \Gamma_i \), \( 1 \leq i \leq k-1 \), of covering transformations of \( P_i \), 
\( 1 \leq i \leq k-1 \), satisfies  
\[
\Gamma_i = \mu_i^{-1} \circ \Gamma_k \circ \mu_i, \quad 1 \leq i \leq k-1,
\]  
the group \( \Gamma \) contains all the groups \( \Gamma_i \), \( 1 \leq i \leq k \), 
 implying as in the proof of Theorem \ref{tt1} that there exist 
a finitely generated subfield $\bm k$ of $\C$ and infinite sets $K_i\subset \P^1({\bm k})$,  $1\leq i \leq k,$ such that the  equality 
\be \l{st} 
(P_k\circ \mu_1)^{-1}(K_1)=\dots =(P_{k}\circ \mu_{k-1})^{-1}(K_{k-1})=P_k^{-1}(K_k)=S \ee
holds.   

Setting \be \l{ifo} \t K_i=P_k^{-1}(K_i), \ \ \ 1\leq i \leq k,\ee we see that \eqref{st} implies the equality 
\be \l{thu} \mu_1^{-1}(\t K_1)=\mu_2^{-1}(\t K_2)=\dots =\mu_{k-1}^{-1}(\t K_{k-1})=\t K_k=S.\ee
On the other hand, since the M\"obius transformations $\mu_i$, $1\leq i \leq k-1,$  belong to $\Gamma$, the set 
$S$ is invariant with respect to these transformations. Thus, \eqref{thu} implies the equality 
$$\t K_1=\t K_2=\dots =\t K_{k-1}=\t K_k=S.$$ It follows now from \eqref{ifo} that 
$$K_1=K_2=\dots =K_{k-1}=K_k=P_k(S).$$ 

The case \( g(R) = 1 \) can be treated using the same approach as in the proof of 
Theorem \ref{tt1}, with appropriate modifications. Namely, we consider the algebraic curve 
\( X \) defined in the proof of Theorem \ref{tt1}, and the group \( \Gamma \subset \Aut(R) \) 
generated by the group \( \Gamma_k \), the automorphisms \( \mu_i \), \( 1 \leq i \leq k-1 \), 
and an arbitrary shift \( \mu \) of infinite order on \( R \). We then define 
\( \widehat{\Gamma} \subset \Aut(X) \) as the group consisting of all elements 
\( \widehat{\sigma} \) with \( \sigma \in \Gamma \), and consider the set \eqref{set}.
 \qed

Notice that, as in the examples given by Theorem \ref{tt1}, the holomorphic maps \( P_i \), \( 1 \leq i \leq k \), provided by Theorem \ref{tt2} may or may not have a nontrivial common compositional right factor, depending on whether the group \eqref{gr1} is trivial.
Furthermore, if the group \eqref{gr2} is finite, then there exist  a holomorphic Galois covering \( A: R \to \mathbb{P}^1(\mathbb{C}) \) and 
a rational function $F$ such that the following holds: \( G \) is the group of covering transformations of $A$, the equality 
\be \l{pop} A= F \circ P_1 = F \circ P_2 = \dots = F \circ P_k \ee holds, and 
$ K = F^{-1}(\h K)$ for some set $\h K$, contained in $\P^1(\bm k)$ for some finitely generated field $\bm k$. 

Indeed, in the considered case, equality \eqref{bn} reduces to the following:
\[
A = F_1 \circ P_k \circ \mu_1 = F_2 \circ P_k \circ \mu_2 = \dots = F_k \circ P_k.
\]
Since \( \mu_i \), \( 1 \leq i \leq k-1 \), belong to \( \Gamma \), we have
\be \l{og} 
A = A \circ \mu_i, \quad 1 \leq i \leq k-1.
\ee
Substituting \( F_i \circ P_k \circ \mu_i \) for \( A \) in the left-hand side and \( F_k \circ P_k \) for \( A \) in the right-hand side of \eqref{og}, we obtain
\[
F_i \circ P_k \circ \mu_i = F_k \circ P_k \circ \mu_i, \quad 1 \leq i \leq k-1,
\]
which implies that
\[
F_1 = F_2 = \dots = F_k.
\]

\end{document}